\documentclass[reqno,12pt]{amsart}

\usepackage{amsmath}
\usepackage{amsfonts}
\usepackage{amssymb}
\usepackage{eufrak}
\usepackage{amscd}
\usepackage{amsthm}
\usepackage{epsfig}
\usepackage{amstext}
\usepackage[all]{xy}

 \newcommand{\K}{{\mathbb K}}

 \newcommand{\Ext}{\operatorname{Ext}}
\newcommand{\Ad}{\operatorname{Ad}}

\newcommand{\id}{\operatorname{id}}

\newcommand{\ev}{\operatorname{ev}}

   \theoremstyle{plain}
   \newtheorem{thm}{Theorem}[section]
   
   \newtheorem{lem}[thm]{Lemma}
   
   \theoremstyle{definition}

   \theoremstyle{remark}

 \numberwithin{equation}{section}

 \title[On the non-invertibility of certain $C^*$-extensions]%
{On the lack of inverses to $C^*$-extensions related to property T groups}
\author{V. Manuilov}
\address{Department of Mechanics and Mathematics, Moscow State
University, Leninskie Gory, Moscow, 119992, Russia}
\email{manuilov@mech.math.msu.su}
\thanks{The first named author acknowledges partial support from
RFFI, grant No. 02-01-00574, and
\mbox{HIII\hspace{-2.3ex}%
\rule{1.9ex}{0.07ex}\,-619.2003.01}}

\author{K. Thomsen}
\address{IMF, Department of Mathematics, Ny Munkegade, 8000 Aarhus C,
Denmark}
\email{matkt@imf.au.dk}
\thanks{}


\begin{document}

 \begin{abstract}
Using ideas of S. Wassermann on non-exact $C^*$-algebras and
property T groups, we show that one of his examples of non-invertible
$C^*$-extensions is not semi-invertible. To prove this, we
show that a certain element vanishes in the asymptotic tensor
product. We also show that a modification of the example gives
a $C^*$-extension which is not even invertible up to homotopy.

 \end{abstract}

\maketitle

\section*{Introduction}

The Brown--Douglas--Fillmore theory of $C^*$-extensions,
\cite{BDF}, works nicely for nuclear $C^*$-algebras because an
extension of a nuclear $C^*$-algebra is always invertible in the
extension semi-group. As a steadily growing number of examples
show, this is not the case for general extensions, cf.
\cite{Anderson},\cite{Kirchberg},\cite{Wass1},\cite{Wass2},\cite{W},\cite{Haagerup},\cite{Ozawa},
etc. In contrast, besides all its other merits, the $E$-theory of
Connes and Higson, \cite{CH}, provides a framework which
incorporates arbitrary extensions of $C^*$-algebras, and in
previous work we have clarified in which way this happens, cf.
\cite{MT2},\cite{MT3}. Specifically, in the $E$-theory setting the
notion of triviality of extensions must we weakened, at least so
far as to consider an extension of $C^*$-algebras
\begin{equation}\label{A1}
\begin{xymatrix}{
0 \ar[r]  &  B \ar[r] &  E \ar[r]^-{q}  &  A  \ar[r]  & 0}
\end{xymatrix}
\end{equation}
to be trivial when it is asymptotically split, by which we mean that there
is an asymptotic homomorphism, \cite{CH},
$\varphi=(\varphi_t)_{t\in[0,\infty)}:A\to E$ such that
$q\circ\varphi_t=\id_A$ for each $t\in[0,\infty)$. When the quotient
$C^*$-algebra $A$ is a suspension, i.e. is of the form
$C_0(\mathbb R) \otimes D$, this is the only change which is needed
to ensure that $E$-theory becomes a complete analogue of the
BDF theory for nuclear $C^*$-algebras. Specifically, when
the quotient $C^*$-algebra is a suspension
and the ideal is stable, every extension is semi-invertible, by which
we mean that it is invertible in the sense corresponding to the
weakened notion of triviality, i.e. one can add an extension to it
so that the result is asymptotically split. Furthermore a
given extension will represent $0$ in
$E$-theory if and only if it can be made asymptotically split by
adding an asymptotically split extension to it. One purpose of
the present paper is to show by example that this nice situation
does not persist when the quotient $C^*$-algebra is not a suspension.
We will show that an extension considered by S.~Wassermann in \cite{W},
and shown by him to be non-invertible, is not
semi-invertible either. By slightly modifying Wassermann's example,
we obtain also an extension which is not even invertible up to homotopy,
giving us the first example of a $C^*$-algebra
for which the semi-group of homotopy classes of extensions by a stable
$C^*$-algebra, in casu the algebra of
compact operators, is not a group. The conclusion
is that the $E$-theory approach to $C^*$-extensions does not completely
save us from the unpleasantness of extensions
without inverses. But unlike the BDF theory, as shown in
\cite{MT3}, in $E$-theory they can be eliminated at the cost of a
single suspension.

The method we adopt in order to show that Wassermann's example from
\cite{W} is not semi-invertible is basically the same as
his, although the verification is somewhat more complicated since
it uses the asymptotic tensor norm, which was introduced
in \cite{MT7}, in place of the minimal tensor norm. To show that a
suitably modified version of the extension is not even invertible
up to homotopy we proceed quite differently in order to bring a
$K$-theoretical obstruction to bear.

\section{Wassermann's extension is not semi-invertible}

\subsection{Wassermann's example}

Let $G$ be an infinite countable discrete group with the property
T of Kazhdan, \cite{Va}. It is a result of Wang, \cite{Wa}, cf.
\cite{W}, that there is at most a countable number of unitary
equivalence classes of finite-dimensional unitary representations
of $G$. As in \cite{W}, we assume that there actually are
infinitely many equivalence classes of such representations.  This
is the case, for example, when $G = SL_3\left(\mathbb Z\right)$.
We fix then a sequence $\pi_k, k = 1,2,3, \dots$, of inequivalent
finite-dimensional irreducible unitary representations of $G$
which contains a representative for each equivalence class of such
representations. Consider the direct sum $\pi =
\oplus_{k=1}^{\infty} \pi_k$ of these representations, acting on
the Hilbert space $H$, and let $B$ be the $C^*$-subalgebra of
${\mathbb L}(H)$ generated by $\left\{\pi(g): g \in G\right\}$.
The $C^*$-subalgebra of ${\mathbb L}(H)$ generated by $B$ and the
ideal $\mathbb K = \mathbb K(H)$ of compact operators on $H$ will
be denoted by $E$. Then $\mathbb K$ is an ideal in $E$ and we
denote the quotient $E/\mathbb K$ by $A$. It was shown in \cite{W}
that the extension
 \begin{equation}\label{C1}
 \begin{xymatrix}{
0 \ar[r]  &  {\mathbb K} \ar[r] &  E \ar[r]^-{q}  &  A  \ar[r]  &
0}
 \end{xymatrix}
 \end{equation}
is not invertible, or not semi-split. We are going to prove that it is
also not semi-invertible.

 \begin{thm}\label{T1}
The extension $($\ref{C1}$)$ is not semi-invertible.
 \end{thm}

We shall elaborate a little
on Wassermann's argument, so let us therefore first outline
it. He shows that (\ref{C1}) does not admit a
completely positive section for the quotient map (i.e. is not
invertible) because the sequence
 \begin{equation}\label{C2}
 \begin{xymatrix}{
  {\mathbb K} \otimes B \ar[r] &  E \otimes_{\text{min}} B \ar[r]  &
A\otimes_{\text{min}} B  }
 \end{xymatrix}
 \end{equation}
is not exact. That (\ref{C2}) is not exact he deduces as follows:
The representation $G \ni g \mapsto \pi(g) \otimes \pi(g)$ of $G$
in $E \otimes_{\text{min}} B$ extends by the universal property of
$C^*(G)$ to a $*$-homomorphism
 $$
 \begin{xymatrix}{
  \Delta : C^*(G) \to  E \otimes_{\text{min}} B .  }
 \end{xymatrix}
 $$
By Theorem 1.10 of \cite{HV}, $G$ is finitely generated.
Let $g_1,g_2, \dots, g_n$ be a set of generators of $G$. We assume that
this set contains the neutral element and is symmetric, i.e. contains
$g_i^{-1}$ for all $i=1,\ldots, n$.
Wassermann shows in the proof of Theorem 6 of
\cite{W}
that there is $\delta>0$ such that the spectrum of
the image of the element
 $$
 \Delta\Bigl(\frac{1}{n} \sum_{i=1}^n g_i\Bigr) \in E\otimes_{\min}B
 $$
in the quotient $E \otimes_{\text{min}} B/ {\mathbb K} \otimes B$
lies in $[-1,1-\delta]\cup\{1\}$ and contains 1, while the
spectrum of its image in $A\otimes_{\text{min}} B$, under the
quotient map of (\ref{C2}), lies in $[-1,1-\delta]$. Thus
(\ref{C2}) is clearly not exact, and it follows that (\ref{C1}) is
not invertible.

In order to adopt this approach in the asymptotic setting it is
crucial that Wassermann's proof of Theorem 6 in \cite{W} gives a tiny
bit of additional information. Recall, \cite{Va},\cite{dHRV}, that
$1$ is an isolated point in the spectrum of
$\frac{1}{n} \sum_{i=1}^n g_i$
in $C^*(G)$, and that the corresponding spectral projection $p$ is
the support projection of the trivial representation. In particular,
$p = h\left(  \frac{1}{n} \sum_{i=1}^n g_i\right)$ for an appropriately
chosen continuous function $h$ on $[-1,1]$, and then it is clear that
Wassermann's argument for Theorem 6 in \cite{W} proves that the image of
$p$
is non-zero in $E \otimes_{\text{min}} B/ {\mathbb K} \otimes B$, but zero
in $A\otimes_{\text{min}} B$. It is this fact, that the non-invertibility
of (\ref{C1}) can be detected by the non-vanishing of a projection in a
certain $C^*$-algebra, which makes it possible to adopt it to the
asymptotic case.

\subsection{Left asymptotic tensor $C^*$-norm}

Let us review the construction of the left asymptotic tensor norm,
cf. \cite{MT7}. Let $\varphi = \left(\varphi_t\right)_{t \in
[1,\infty)} : A \to {\mathbb L}(H_1)$ be an asymptotic
homomorphism from $A$ to the $C^*$-algebra of bounded operators
${\mathbb L}(H_1)$ of some Hilbert space $H_1$, in the following
referred to as \emph{an asymptotic representation} of $A$, and let
$\pi : B \to {\mathbb L}(H_2)$ be a (genuine) representation of
$B$. $\varphi$ and $\pi$ define in the natural way two commuting
$*$-homomorphisms,
 $$
A \to C_b\left([ 1, \infty), {\mathbb L}(H_1 \otimes H_2)\right)/
C_0\left([ 1, \infty), {\mathbb L}(H_1 \otimes H_2)\right)
 $$
and
 $$
B  \to C_b\left([ 1, \infty), {\mathbb L}(H_1 \otimes H_2)\right)/
C_0\left([ 1, \infty), {\mathbb L}(H_1 \otimes H_2)\right),
 $$
giving rise to a $*$-homomorphism
 $$
\varphi \odot \pi : A \odot B \to  C_b\left([ 1, \infty), {\mathbb
L}(H_1 \otimes H_2)\right)/ C_0\left([ 1, \infty), {\mathbb L}(H_1
\otimes H_2)\right).
 $$
The left asymptotic tensor norm $\|\cdot\|_{\lambda}$, defined on
$A\odot B$, is
 $$
\|c\|_{\lambda} = \sup_{\varphi, \pi} \left\|\varphi \odot \pi(c)\right\|,
 $$
where the supremum is taken over all asymptotic representations of $A$
and all representations of $B$. On a linear combination,
$c =\sum_{i=1}^m a_i \otimes b_i$, of simple tensors,
 $$
\|c\|_{\lambda} = \sup_{\varphi, \pi} \Bigl( \limsup_{t \to \infty}
\Big\|\sum_{i=1}^m \varphi_t(a_i) \otimes \pi(b_i)\Big\|\Bigr) .
 $$
Let $A \otimes_{\lambda} B$ be the completion of $A \odot B$
with respect to the norm $\|\cdot\|_\lambda$.

It is a convenient fact that the left asymptotic tensor norm can be
calculated using only a single representation of $B$:

 \begin{lem}\label{C7}
Let $\pi' : B \to {\mathbb L}(H')$ be a faithful representation of
$B$ such that $\pi'(B) \cap \mathbb K(H') = \{0\}$. Then
 $$
\|c\|_{\lambda} = \sup_{\varphi} \left\|\varphi \odot \pi'(c)\right\|
 $$
for all $c \in A \odot B$.
 \end{lem}

 \begin{proof} Let $\rho: B \to {\mathbb L}(H_2)$ be an arbitrary representation
of $B$. We claim that
 \begin{equation}\label{C9}
\left\|\varphi \odot \rho( c)\right\| \leq \left\|\varphi \odot
\pi'( c)\right\|
 \end{equation}
for every asymptotic representation $\varphi$ of $A$ and every $c
\in A \odot B$. To prove this, let $\varepsilon > 0$ and write $c
=\sum_{i=1}^m a_i \otimes b_i$. By Voiculescu's non-commutative
Weyl-von Neumann theorem, \cite{V}, there is an isometry $V : H_2
\to H'$ such that $\sum_{i=1}^m \left\|a_i\right\|\left\|\rho(b_i)
-
  V^*\pi'(b_i)V\right\|
\leq \varepsilon$. Since $\limsup_{t \to \infty}
\left\|\varphi_t(a_i)\right\| \leq \|a_i\|$, it follows that
 \begin{eqnarray*}
\left\|\varphi \odot \rho( c)\right\| &\leq &\limsup_{t \to
\infty} \Big\| \sum_{i=1}^m \varphi_t(a_i) \otimes V^*
\pi'(b_i)V\Big\| + \varepsilon \\ &\leq & \limsup_{t \to
\infty}\Big\| (1 \otimes V^*) \Big(\sum_{i=1}^m \varphi_t(a_i)
\otimes \pi'(b_i)\Big)(1 \otimes V^*) \Big\| + \varepsilon \\
&\leq & \limsup_{t \to \infty} \Big\| \sum_{i=1}^m \varphi_t(a_i)
\otimes \pi'(b_i) \Big\| + \varepsilon,
 \end{eqnarray*}
 proving (\ref{C9}) and hence the lemma.
 \end{proof}

We show now how the asymptotic tensor norm
can be used in proving non-semi-invertibility of an extension.

Note that, thanks to the exactness of the maximal tensor product,
$E\otimes_{\min}B/\K\otimes B$ is a quotient of
$A\otimes_{\max}B$. On the other hand, $A\otimes_{\min}B$ is the
quotient of $E\otimes_{\min}B/\K\otimes B$. Therefore $A\odot B$
is a dense subspace in $E\otimes_{\min}B/\K\otimes B$. We denote
the norm on $A\odot B$ inherited from $E\otimes_{\min}B/\K\otimes
B$ by $\|\cdot\|_E$. Since this norm is a cross-norm, we may write
view $E\otimes_{\min}B/\K\otimes B$ as a tensor product of $A$ and $B$ and
write $E\otimes_{\min}B/\K\otimes B=A\otimes_E B$.
Recall that $A\odot B$ is dense in
$A\otimes_\lambda B$ as well.

 \begin{lem}\label{exact_seq}
Suppose that there exists $c\in A\odot B$ such that
$\|c\|_E>\|c\|_\lambda$. Then the extension $($\ref{C1}$)$ is not
semi-invertible.
 \end{lem}

 \begin{proof}
The idea of the proof is borrowed from \cite{Wass2}.
Suppose the contrary, i.e. that (\ref{C1}) is semi-invertible. Then
there exists an extension
 $
\begin{xymatrix}{
0\ar[r]& \K\ar[r]& E'\ar[r]^-{q'}& A\ar[r]& 0 }
\end{xymatrix}
 $
and an asymptotic splitting $s=(s_t)_{t\in{[0,\infty)}}:A\to C$,
where $D\subset M_2({\mathbb L}(H))$ is the $C^*$-subalgebra
 $$
D=\left\lbrace\left(\begin{array}{cc}e&b_1\\b_2&e'\end{array}\right):
b_1,b_2\in \K, e\in E,e'\in E', q(e)=q'(e')\right\rbrace.
 $$
By definitition of the left asymptotic tensor norm there is an asymptotic
homomorphism
$s_t\otimes_\lambda \id_B:A\otimes_\lambda B\to D\otimes_{\min}
B$ with the property that
$$
\lim_{t \to \infty} s_t\otimes_\lambda \id_B\left( a \odot b\right) -
s_t(a) \odot b = 0
$$
on simple tensors. Let $d:D\to E$ be the completely positive contraction
given by
$d\left(\begin{array}{cc}e&b_1\\b_2&e'\end{array}\right)=e$. Then
the map $d\otimes\id_B:D\otimes_{\min}B\to E\otimes_{\min}B$ is a
well-defined contraction. Let $q_B:E\otimes_{\min}B\to
E\otimes_{\min}B/\K\otimes B=A\otimes_E B$ be the quotient map.
Consider the composition
 \begin{equation*}\label{D1}
r_t=q_B\circ (d\otimes\id_B)\circ
(s_t\otimes_\lambda\id_B):A\otimes_\lambda B\to
A\otimes_E B.
 \end{equation*}
The maps $q_B$ and $d\otimes\id_B$ are contractions and the family
$(s_t\otimes_\lambda\id_B)_{t \in [1,\infty)}$ is asymptotically
contractive, so the family $(r_t)_{t \in [1,\infty)}$ is asymptotically
contractive. Since
$\lim_{t \to \infty} r_t(c) -c = 0$, it follows that
$\|c\|_E=\lim\sup_{t\to\infty}\|r_t(c)\|_E\leq \|c\|_\lambda$.
The contradiction to $\|c\|_E>\|c\|_\lambda$ completes the proof.
 \end{proof}

Let $f(t)$ be a polynomial $f(t) = \frac{1}{4}t^2 + \frac{1}{2}t +
\frac{1}{4}$, and set $x=f\left(\frac{1}{n} \sum_{i=1}^ng_i\right)\in
C^*(G)$. Then $0 \leq x \leq 1$ and $1$ is an isolated point in the
spectrum of $x$.
Put $\Delta(x)\in A\odot B$. As pointed out above, Wassermann has shown
 that the spectrum of the element
$\Delta(\frac{1}{n}\sum_{i=1}^ng_i)\in A\odot B$ in the
quotient $E\otimes_{\min}B/\K\otimes B$ contains 1, and it follows that
$\|\Delta(x)\|_E=1$. By Lemma \ref{exact_seq}, Theorem \ref{T1} will
follow if we show that
 \begin{equation}\label{D2}
\|\Delta(x)\|_\lambda <1.
 \end{equation}
Let $\Delta_{\lambda} : C^*(G) \to A \otimes_{\lambda} B$ be the
$*$-homomorphism determined by the condition that
$\Delta_{\lambda}(g) = q\left(\pi(g)\right) \otimes \pi(g), g \in G$.
The desired conclusion, (\ref{D2}),
is then equivalent to
 $$
\left\|\Delta_{\lambda}\left(p\right)\right\|_{\lambda}  = 0,
 $$
because $1$ is isolated in the spectrum of $x$.

\subsection{Calculation of $\|\Delta_\lambda(p)\|_\lambda$}

 \begin{lem}\label{norm}
One has $\|\Delta_\lambda(p)\|_\lambda=0$.
 \end{lem}

 \begin{proof}
Set $H' = \oplus_{i=1}^{\infty} H$ and let $i_{\infty} : B \to
{\mathbb L}\left(H'\right)$ be the infinite sum of copies of the
inclusion $B \subseteq {\mathbb L}(H)$. Then
 \begin{equation*}
\left\|c\right\|_{\lambda} =  \sup_{\varphi}
\left\|\varphi \odot i_{\infty}( c)\right\|
 \end{equation*}
for all $c \in A \odot B$ by Lemma \ref{C7}. Let $\varepsilon \in
(0,\frac{1}{100})$. There is then an asymptotic representation
$\varphi : A \to {\mathbb L}(H_1)$ and an equi-continuous
asymptotic representation $\Phi: A \otimes_{\lambda} B \to
{\mathbb L}(H_1 \otimes H')$ such that
 \begin{equation}\label{C10}
\limsup_{t \to \infty} \left\|\Phi_t\left(\Delta_{\lambda}(p)
\right)\right\| \geq \|\Delta_{\lambda}(p)\|_{\lambda}
-\varepsilon,
 \end{equation}
and
 \begin{equation}\label{C11}
\lim_{t \to \infty} \Big\|\Phi_t\left(c\right) - \sum_{i=1}^m
\varphi_t(a_i) \otimes i_{\infty}(b_i)\Big\| = 0
 \end{equation}
for all $c = \sum_{i=1}^m a_i \otimes b_i \in A \odot B$. For each
$k$, let $q_k'$ be the orthogonal projection onto the support in
$H$ of the representation $\pi_k$, and let $q_k =
\prod_{i=1}^{\infty} q_k' \in {\mathbb L}(H')$ be the infinite
repeat of $q_k'$. Note that each $1_{H_1} \otimes q_k$ commutes
with $\sum_{i=1}^m \varphi_t(a_i) \otimes i_{\infty}(b_i)$ for all
$t$ and all $c = \sum_{i=1}^m a_i \otimes b_i \in A \odot B$. By
approximating $\Delta_{\lambda}(p)$ with elements from $A \odot B$
we can find an element $z =  \sum_{i=1}^m a_i \otimes b_i \in A
\odot B$ such that
 \begin{equation}\label{C13}
\limsup_{t \to \infty} \Big\|
\Phi_t\left(\Delta_{\lambda}(p)\right) - \sum_{i=1}^m
\varphi_t(a_i) \otimes i_{\infty}(b_i)\Big\| < \varepsilon .
 \end{equation}
To simplify notation, set $z_t = \sum_{i=1}^m \varphi_t(a_i) \otimes
i_{\infty}(b_i)$, and $y_t = \frac{1}{2}\left(z_t + z_t^*\right)$.
Since $\Phi$ is an asymptotic homomorphism and $\Delta_{\lambda}(p)$
a projection it follows from (\ref{C13}) that for some $T > 0$,
 \begin{equation}\label{C14}
\left\|y_t^2 - y_t\right\| \leq 5 \varepsilon
 \end{equation}
when $t \geq T$. It follows that
 $$
\big\|\left(\left(1_{H_1} \otimes q_k\right) y_t\right)^2 -
\left(1_{H_1} \otimes q_k\right)y_t\big\| \leq 5 \varepsilon
 $$
for all $t > T$. Since $5 \varepsilon < 1/4$, we find that the
characteristic function $h = 1_{[1/2,\infty)}$ is continuous on
the spectrum of $y_t$ and on the spectrum of each $\left(1_{H_1} \otimes
q_k\right) y_t$ when $t > T$. It follows that $h(y_t)$ and
$h\left(\left(1_{H_1} \otimes q_k\right) y_t\right)$ are
projections for all $k$ and all $t > T$. We claim that
 \begin{equation}\label{C16}
h(y_t) = 0
 \end{equation}
for all $t > T$. If not there is some $t_0 > T$ such that $h(y_{t_0})
\neq 0$. There must then be a $k$, which we now fix, such that
$h\left(\left(1_{H_1} \otimes q_k\right) y_{t_0}\right) \neq 0$
since $\sum_{i} 1_{H_1} \otimes q_i = 1$. But then
 \begin{equation}\label{C17}
\left\|h\left(\left(1_{H_1} \otimes q_k\right)y_t\right)\right\| =1
 \end{equation}
for all $t > T$ since $h\left(\left(1_{H_1} \otimes
q_k\right)y_t\right)$ varies norm-continuously with $t$ and is a
projection for all $t > T$. Let $\rho_k : B \to
C^*\left(\pi_n(G)\right)$ be the finite-dimensional representation
of $B$ obtained by restricting the elements of $B$ to the subspace
of $H$ supporting the representation $\pi_k$ of $G$. There is then
a representation $\mu : C^*\left(\pi_n(G)\right) \to {\mathbb
L}(H')$ such that
 \begin{equation}\label{C19}
\mu \circ \rho_k(b) = q_k i_{\infty}(b)
 \end{equation}
for all $b \in B$. Furthermore, there is an equi-continuous
asymptotic homomorphism $\psi : A \otimes C^*\left(\pi_n(G)\right)
\to {\mathbb L}(H_1 \otimes H')$ such that
 \begin{equation}\label{C18}
\lim_{t \to \infty} \Big\|\psi_t(c) - \sum_{i=1}^m \varphi_t(a_i)
\otimes \mu(x_i)\Big\| = 0
 \end{equation}
for all $c = \sum_{i=1}^m a_i \otimes x_i \in A \odot
C^*\left(\pi_n(G)\right)$. Note that $\id_A \otimes \rho_k :
A \odot B \to A \odot C^*\left(\pi_n(G)\right)$ extends to a
$*$-homomorphism $\kappa : A \otimes_{\lambda} B \to  A \otimes
C^*\left(\pi_n(G)\right)$. It follows from (\ref{C11}), (\ref{C19})
and (\ref{C18}) that
 \begin{equation}\label{C20}
\lim_{t \to \infty} \left\|\psi_t \circ \kappa(d) - \left(1_{H_1}
\otimes q_k\right) \Phi_t(d) \right\| = 0
 \end{equation}
for all $d \in A \otimes_{\lambda} B$. Since $\kappa$ factors
through $A \otimes_{\min} B$, we know from \cite{W} that
$\left\|\kappa\left(\left(\Delta_{\lambda}(p)\right)\right)\right\|  =
0$.
It follows therefore from (\ref{C20}) that
 $$
\limsup_{t \to \infty} \left\|\left(1_{H_1} \otimes q_k\right)
\Phi_t\left(\left( \Delta_{\lambda}(p) \right) \right) \right\| = 0,
 $$
and then by use of (\ref{C13}) that
 $$
\limsup_{t \to \infty} \left\| \left(1_{H_1} \otimes q_k\right)
y_t \right\| \leq \varepsilon .
 $$
Since $\varepsilon < 1/2$, this contradicts (\ref{C17}), and we
conclude that (\ref{C16}) must hold. Combined with (\ref{C14}) we
find that the spectrum of $y_t$ is contained in $[- 1/2,1/2]$, and
hence that $\left\|y_t\right\| \leq 1/2$. It follows then from
(\ref{C13}) that
 $$
\limsup_{t \to \infty} \left\|
\Phi_t\left(\Delta_{\lambda}(p)\right)\right\| \leq 1/2 +
\varepsilon  < 1 .
 $$
Since $\Delta_{\lambda}(p)$ is a projection, we deduce first that
$\limsup_{t \to \infty} \left\|
\Phi_t\left(\Delta_{\lambda}(p)\right)\right\| = 0$, and then from
(\ref{C10}) that $\left\|\Delta_{\lambda}(p)\right\|_{\lambda} =
0$.
 \end{proof}

\subsection{Some remarks}

Theorem \ref{T1} means that it is not possible to add an extension
of $A$ by $\mathbb K$ to (\ref{C1}) such that the resulting
extension admits an asymptotic homomorphism consisting of sections
for the quotient map. In particular, the extension
(\ref{C1}) itself does not admit such a family of sections; a
fact, which may seem slightly surprising because the extension is
clearly quasi-diagonal and there is an obvious sequence $s_n : A
\to E$, $n = 1,2, \dots$, of maps, each of which is a section for
the quotient map such that they form a \emph{discrete} asymptotic
homomorphism. It was therefore no coincidence that the
connectedness of the parameter space, $[0,\infty)$, was used at a
crucial point in the proof above.

In \cite{MT7} we raised the question, if the left asymptotic tensor
product is associative. It follows from Lemma \ref{norm} that the
answer is negative.



We have looked through all known (to us) examples of
non-invertible extensions to check if they are semi-invertible or
not. The examples of Kirchberg, \cite{Kirchberg}, are
semi-invertible by results of \cite{MT3}. Another example of
Wassermann, \cite{Wass1}, can be shown to not to be semi-invertible
by the same method as here. Unfortunately, we know nothing about
semi-invertibility of other examples.

\section{Homotopy non-invertibility}

\subsection{A modification of the Wassermann's extension}

To give an example of an extension which is not only not
semi-invertible, but also not even homotopy invertible, we modify
the extension (\ref{C1}) as follows. Let $d_i$ be the dimension of
the Hilbert space $H_i$ on which the representation $\pi_i$ acts.
Let $n_i$ be a sequence of integers such
that $\lim_{i \to \infty} \frac{n_i}{d_i} = \infty$. For each $i$
we let $n_i \cdot \pi_i$ be the direct sum of $n_i$ copies of the
representation $\pi_i$, and let
 $$
\pi'= \oplus_{i=1}^{\infty} n_i \cdot \pi_i
 $$
be the direct sum of the resulting sequence of representations,
acting on the Hilbert space $H$. Let $E'$ be the $C^*$-subalgebra
of $\mathbb L(H)$ generated by $\{\pi'(g): g \in G\}$ and by
$\mathbb K$, the compact operators on $H$. Set $A' = E'/\mathbb
K$.

 \begin{thm}\label{homtop}
The extension
 \begin{equation}\label{hom3}
\begin{xymatrix}{
0 \ar[r]  &  {\mathbb K} \ar[r] &  E' \ar[r]  &  A'  \ar[r]  & 0}
\end{xymatrix}
 \end{equation}
is not invertible in $\Ext_h(A',\mathbb K)$.
 \end{thm}

\begin{proof} To show that (\ref{hom3}) is not invertible up to homotopy,
let $\varphi : A' \to Q(\mathbb K)$ be the Busby invariant of
(\ref{hom3}), and assume to reach a contradiction that $\psi : A'
\to Q(\mathbb K)$ is an extension such that $\varphi \oplus \psi$
is homotopic to $0$. Let $V_1,V_2$ be isometries in $\mathbb L(H)$
such that $V_1V_1^* + V_2V_2^* = 1$, and set $\lambda (a) = \Ad
q(V_1) \circ \varphi(a) + \Ad q(V_2) \circ \psi(a)$, $a \in A'$,
where $q:{\mathbb L}(H)\to Q(\K)$ is the quotient map. There is
then a commuting diagram
\begin{equation}\label{hom}
\begin{xymatrix}{
0 \ar[r]  &  {\mathbb K} \ar[r] &  {\mathcal E} \ar[r]  &  A'
\ar[r]  & 0\\ 0 \ar[r]  & { I \mathbb K} \ar[u]^-{\ev_1}
\ar[d]_-{\ev_0} \ar[r]  & {\mathcal E'} \ar[u] \ar[d] \ar[r] &  A'
\ar[r] \ar@{=}[u] \ar@{=}[d] & 0\\ 0 \ar[r]  &  {\mathbb K} \ar[r]
&  {\mathbb K} \oplus A'  \ar[r]  &  A'  \ar[r]  & 0, }
\end{xymatrix}
\end{equation}
where $\lambda$ is a the Busby invariant of the upper extension,
$I\mathbb K = C[0,1] \otimes \mathbb K$ and $\ev_s : I\mathbb K
\to \mathbb K$ is evaluation at $s \in [0,1]$. Set $D =
\prod_{k=1}^{\infty}
\mathbb L(H_k)$. By tensoring with $D$ we obtain from (\ref{hom})
the commuting diagram
\begin{equation}\label{hom1}
\begin{xymatrix}{
0 \ar[r]  &  {\mathbb K}\otimes D \ar[r] &  {\mathcal E}
\otimes_{\min} D \ar[r]  &  {\mathcal E} \otimes_{\min} D/
{\mathbb K}\otimes D    \ar[r]  & 0\\ 0 \ar[r]  & { I \mathbb K}
\otimes D \ar[u]^-{\ev_1}   \ar[d]_-{\ev_0} \ar[r]  & {\mathcal
E'}\otimes_{\min} D  \ar[u] \ar[d] \ar[r] &  {\mathcal
E'}\otimes_{\min} D/ {{ I \mathbb K} \otimes D} \ar[r]
\ar[u]^-{p_0} \ar[d]_-{p_1} & 0\\ 0 \ar[r]  &  {\mathbb K}\otimes
D \ar[r] &  ({\mathbb K}\otimes D) \oplus (A \otimes_{\min} D)
\ar[r]  &   A' \otimes_{\min} D   \ar[r]  & 0. }
\end{xymatrix}
\end{equation}

Let $\overline{\ev}_s : M(I\mathbb K) \to M(\mathbb K)$ and
$\widehat{\ev}_s : Q(I\mathbb K) \to Q(\mathbb K)$ be the
$*$-homomorphisms induced by $\ev_s$. Denote $ \mathcal E
\otimes_{\min} D/ {\mathbb K}\otimes D $ and $\mathcal
E'\otimes_{\min} D/ {{ I \mathbb K} \otimes D}$ by $A'
\otimes_{\mathcal E} D$ and $A' \otimes_{\mathcal E'} D$,
respectively. The Busby invariant of the middle extension of
(\ref{hom1}) is a $*$-homomorphism $\varphi' : A'
\otimes_{\mathcal E'} D \to Q(I \mathbb K \otimes D)$ such that
$\widehat{\ev}_1 \circ \varphi' = \mu \circ p_0$, where $\mu : A'
\otimes_{\mathcal E} D \to Q(I\mathbb K \otimes D)$ is the Busby
invariant of the upper extension in (\ref{hom1}), while
$\widehat{\ev}_0 \circ \varphi' = 0$.

By the Bartle-Graves selection theorem there are continuous
sections $\chi : Q(I\mathbb K \otimes D) \to M(I\mathbb K \otimes
D)$ and $\chi_k : Q(I\mathbb K \otimes D_k) \to M(I\mathbb K
\otimes D_k)$ for the quotient maps $ M(I\mathbb K \otimes D) \to
Q(I\mathbb K \otimes D)$ and $ M(I\mathbb K \otimes D_k) \to
Q(I\mathbb K \otimes D_k)$, respectively, for all $ k = 1,2,3,
\dots$. We can choose these maps to be self-adjoint and such that
$\|\chi(x)\| \leq 2\|x\|$, $x \in Q(I\mathbb K \otimes D)$, and
$\|\chi_k(y)\| \leq 2\|y\|$, $y \in Q(I\mathbb K \otimes D_k)$,
for all $k$. Set $D_k = \mathbb L(H_k)$, and let $p_k :
D=\prod_{k=1}^\infty D_k \to D_k$ be the canonical projection. The
map $\id_A \otimes p_k: A' \odot D \to A' \odot D_k$ extends to a
$*$-homomorphism $\id_A \otimes p_k : A' \otimes_{\mathcal E'} D
\to A' \otimes D_k$. Let $\overline{\id_{I \mathbb K} \otimes p_k}
: M(I\mathbb K \otimes D) \to M(I\mathbb K \otimes D_k)$ be the
unique $*$-homomorphism extending $\id_{I \mathbb K} \otimes p_k :
I \mathbb K \otimes D \to I \mathbb K \otimes D_k$, and
$\widehat{\id_{I \mathbb K} \otimes p_k} : Q(I\mathbb K \otimes D)
\to Q(I\mathbb K \otimes D_k)$ the resulting $*$-homomorphism. Let
$\Phi : A' \to Q(I\mathbb K)$ be the Busby invariant of the middle
extension of (\ref{hom}).

We denote by $\Phi \hat{\otimes} \id_{D_k}$ the $*$-homomorphism
$A' \otimes D_k \to Q(I\mathbb K \otimes D_k)$ obtained by
composing $\Phi \otimes \id_{D_k} : A' \otimes D_k \to Q(I\mathbb
K)\otimes D_k$ with the canonical embedding $Q(I\mathbb K)\otimes
D_k \subseteq Q(I\mathbb K \otimes D_k)$. By checking on simple
tensors one finds that
 $$
\widehat{\id_{I \mathbb K} \otimes p_k} \circ \varphi' = \left(
\Phi \hat{\otimes} \id_{D_k}\right) \circ \left( \id_{A'} \otimes
p_k\right),
 $$
which implies that
\begin{equation}\label{urt}
\left(\overline{\id_{I \mathbb K} \otimes p_k}\right) \circ \chi
\circ \varphi'(x) - \chi_k \circ \left(\Phi
\hat{\otimes}\id_{D_k}\right) \circ ( \id_{A'} \otimes p_k)(x) \in
I\mathbb K \otimes D_k
\end{equation}
for all $k$ and all $x \in A \otimes_{\mathcal E'} D$.
Let $\overline{\pi_i}$ be
the representation of $G$ contragredient to $\pi_i$. The
representation $g \mapsto q\left(\pi'(g)\right) \otimes \left(
\prod_{i=1}^{\infty} \overline{\pi}_i\right)(g)$ of $G$ into $A'
\odot D$ gives rise to a $*$-homomorphism $\Delta' : C^*(G) \to A'
\otimes_{\mathcal E'} D$. Set $Q = \Delta'(p)$, where $p$, as
above, is the spectral projection of the element $\frac{1}{n}
\sum_{i=1}^n g_i$ corresponding to the set $\{1\}$. Since $\pi_i$
is inequivalent to $\pi_k$ for all except finitely many $i$, it follows from
\cite{W}, Lemma 1, that $ \id_{A'} \otimes p_k(Q) = 0$ for all
$k$. It follows therefore from (\ref{urt}) that
\begin{equation}\label{urt2}
\left(\overline{\id_{I \mathbb K} \otimes p_k}\right) \circ \chi
\circ \varphi'(Q) \in I\mathbb K \otimes D_k
\end{equation}
 for all $k$. Let $P_I$ and $P$ denote the $*$-homomorphisms
$P_I = \prod_{i=1}^{\infty}\overline{\id_{I \mathbb K} \otimes
p_k} : M\left(I\mathbb K \otimes D\right) \to \prod_{i=1}^{\infty}
M\left(I\mathbb K \otimes D_k\right)$ and $P =
\prod_{i=1}^{\infty}\overline{\id_{\mathbb K} \otimes p_k} :
M\left(\mathbb K \otimes D\right) \to \prod_{i=1}^{\infty}
M\left(\mathbb K \otimes D_k\right)$, respectively. Put
 $$
N_I=P_I^{-1}\Bigl(\prod\nolimits_{i=1}^{\infty} {I\mathbb K}
\otimes D_k\Bigr)\subset M({I \mathbb K} \otimes D)
 $$
and
 $$
N=P^{-1}\Bigl(\prod\nolimits_{i=1}^{\infty} {\mathbb K} \otimes
D_k\Bigr)\subset M({\mathbb K} \otimes D).
 $$
It follows from (\ref{urt2}) that
 $$
\chi \circ \varphi'(Q) \in N_I.
 $$
Note that $I\mathbb K \otimes D$ is an ideal in $N_I$ and $\mathbb
K \otimes D$ is an ideal in $N$. We denote the quotients
$N_I/I\mathbb K \otimes D$ and $N/\mathbb K \otimes D$ by $R_I$
and $R$, respectively.
 Note that $R_I \subseteq Q\left(I \mathbb K \otimes
D\right)$ and that $\varphi'(Q) \in R_I$. Evaluation at $s \in
[0,1]$ induces a $*$-homomorphism $E_s : {M\left({I \mathbb K}
\otimes D\right)} \to {M\left({ \mathbb K} \otimes D\right)}$ with
the property that $E_s(N_I)=N$, so we get a $*$-homomorphism
$\widehat{E}_s : R_I \to R$ induced by $E_s$ for each $s \in
[0,1]$. To proceed with the proof, we need some calculations in
$K$-theory.

\subsection{$K$-theory calculations}

Consider the extensions
\begin{equation}\label{urt11}
\begin{xymatrix}{
0 \ar[r] &  {I\mathbb K} \otimes D   \ar[r]  & N_I \ar[r] &  R_I
\ar[r] & 0 }
\end{xymatrix}
\end{equation}
and
\begin{equation}\label{urt11a}
\begin{xymatrix}{
0 \ar[r] &  {\mathbb K} \otimes D   \ar[r]  & N \ar[r] &  R \ar[r]
& 0 }.
\end{xymatrix}
\end{equation}
 The map $\prod_{i=1}^{\infty} {p_i}_* : K_0\left(I\mathbb K
\otimes D\right) \to \prod_{i=1}^{\infty} K_0\left( I\mathbb K
\otimes D_i\right)$ is injective by Lemma 3.2 of \cite{DE} and the
map $\prod_{i=1}^{\infty} {p_i}_* : K_1\left(I\mathbb K \otimes
D\right) \to \prod_{i=1}^{\infty} K_1\left( I\mathbb K \otimes
D_i\right)$ is injective by Lemma 3.3 of \cite{DE}. In particular,
$K_1\left( {I\mathbb K} \otimes D\right) = 0$. Therefore the
extension (\ref{urt11}) gives us a commuting diagram
 $$
\begin{xymatrix}{
0 \ar[d]  &  0 \ar[d] \\
 K_0\left({I\mathbb K} \otimes D\right)  \ar[d] \ar@{=}[r] &
K_0\left({I\mathbb K} \otimes D\right) \ar[d]  \\
 K_0\left(N_I\right)  \ar[r]^-{{P_I}_*} \ar[d] &  {\prod_{i=1}^{\infty}
K_0\left({I\mathbb K}
 \otimes D_k\right)} \ar[d] \\
K_0\left(R_I\right)  \ar[r]  \ar[d]  &  {\bigl(\prod_{i=1}^{\infty}
K_0\left({I\mathbb K} \otimes D_k\right)\bigr)/K_0\left({I\mathbb K}
\otimes D\right) } \ar[d] \\ 0 & 0}
\end{xymatrix}
 $$
 Throwing away the interval, we get also a commuting diagram
 $$
\begin{xymatrix}{
0 \ar[d]  &  0 \ar[d] \\
 K_0\left({\mathbb K} \otimes D\right)  \ar[d] \ar@{=}[r] &
K_0\left({\mathbb K} \otimes D\right) \ar[d]  \\
 K_0\left(N\right)  \ar[r]^-{P_*} \ar[d] &  {\prod_{i=1}^{\infty}
K_0\left({\mathbb K}
 \otimes D_k\right)} \ar[d] \\
K_0\left(R\right)  \ar[r]  \ar[d]  &  {\bigl(\prod_{i=1}^{\infty}
K_0\left({\mathbb K} \otimes D_k\right)\bigr)/K_0\left({\mathbb K}
\otimes D\right) } \ar[d] \\ 0 & 0}
\end{xymatrix}
 $$
Evaluation at any $s \in [0,1]$ induces an isomorphism
 $$
{\left(\prod\nolimits_{i=1}^{\infty} K_0\left({I\mathbb K} \otimes
D_k\right)\right)/K_0\left({I\mathbb K} \otimes D\right) } \to
{\left(\prod\nolimits_{i=1}^{\infty} K_0\left({\mathbb K} \otimes
D_k\right)\right)/K_0\left({\mathbb K} \otimes D\right) }
 $$
in the obvious way and the diagram
 $$
\begin{xymatrix}{
K_0\left(R_I\right) \ar[r]^-{{\widehat{E}_s}{}_*} \ar[d] & K_0(R)
\ar[d] \\
 {\bigl(\prod_{i=1}^{\infty} K_0\left({I\mathbb K} \otimes
D_k\right)\bigr)/K_0\left({I\mathbb K}
 \otimes D\right) } \ar[r]  &  {\bigl(\prod_{i=1}^{\infty}
 K_0\left({\mathbb K} \otimes D_k\right)\bigr)/K_0\left({\mathbb K}
\otimes D\right) } }
\end{xymatrix}
 $$
commutes for every $s \in [0,1]$. Let $x$ be the image in
${\bigl(\prod_{i=1}^{\infty} K_0\left({\mathbb K} \otimes
D_k\right)\bigr)/K_0\left({\mathbb K} \otimes D\right)}$ of the
element $[\varphi'(Q)] \in K_0\left(R_I\right)$. Since
$\widehat{\ev}_0 \circ \varphi' = 0$, we conclude that
$\widehat{E}_0 \left( \varphi'(Q)\right) = 0$, which leads to the
conclusion that
\begin{equation}\label{urt18}
x = 0.
\end{equation}
As we shall see, we get a different result when we consider the
case $s =1$. Let $r : N\to R$ be the quotient map. Set $W_i = V_i
\otimes 1_D \in M(\mathbb K \otimes D)$, $i = 1,2$. Then
\begin{equation}\label{urt18a}
\widehat{E}_1 \circ \varphi'(Q) = r\left( W_1eW_1^* +
W_2aW_2^*\right),
\end{equation}
where $e$ is the spectral projection of $\frac{1}{n} \sum_{i=1}^n
\pi'(g_i) \otimes \left( \oplus_{k=1}^{\infty}
\overline{\pi}_k\right)(g_i)$ corresponding to $\{1\}$, and $a
\geq 0$ is some lift in $N$ of a projection in $R\subset Q\left(
\mathbb K \otimes D\right)$. Since $W_iN\subset N$ and
$W_i^*N \subset N$, the $W_i$'s define multipliers, first of $N$,
and then of $R$. It follows therefore from (\ref{urt18a}) that
 $$
\left[\widehat{E}_1 \circ \varphi'(Q)\right] = [r(e)] + [r(a)]
 $$
in $K_0(R)$.

Consider the extension
\begin{equation}\label{urt19}
\begin{xymatrix}{
0 \ar[r] & {\mathbb K \otimes D} \ar[r]  & N^+ \ar[r]^-{r^+} & R^+
\ar[r]  &  0}
\end{xymatrix}
\end{equation}
obtained from the extension (\ref{urt11a}) by unitalizing. It
follows from Lemma
9.6 of \cite{E} that there are natural numbers $n,m$ such that
$r(a) \oplus 1_n \oplus 0_m \in M_{1+n+m}\left(R^+\right)$ can be
lifted to a projection $f_1 \in M_{1+n+m}\left(N^+\right)$. Note
that the image of $f_1$ in $ M_{1+n+m}\left(\mathbb C\right)$
under the canonical surjection $ M_{1+n+m}\left(N^+\right) \to
M_{1+n+m}\left(\mathbb C\right)$ is a projection of rank $n$. (We
use here that the proof in \cite{E} works equally well when the
assumption that the ideal is AF is replaced by the weaker
assumption, valid in (\ref{urt19}), that the ideal has trivial
$K_1$-group.) There is a commuting diagram
 $$
\begin{xymatrix}{
0 \ar[d]  &  0 \ar[d] \\
 {K_0\left(\mathbb K \otimes D\right)}  \ar[d] \ar@{=}[r] &
{K_0\left(\mathbb K \otimes D\right)} \ar[d]  \\
 {K_0\left(N^+\right)}  \ar[r]^-{P^+_*} \ar[d] &  {\prod_{i=1}^{\infty}
 K_0\left(\left(\mathbb K \otimes D_k\right)^+\right)} \ar[d] \\
K_0\left(R^+\right)  \ar[r]  \ar[d]  &  {\bigl(\prod_{i=1}^{\infty}
K_0\left(\left(\mathbb K \otimes
D_k\right)^+\right)\bigr)/K_0\left(\mathbb K \otimes D\right) } \ar[d]
\\ 0 & 0}
\end{xymatrix}
 $$
with exact columns. Thus the image of $x$ in
${\bigl(\prod_{i=1}^{\infty} K_0\left(\left(\mathbb K \otimes
D_k\right)^+\right)\bigr)/K_0\left(\mathbb K \otimes D\right) }$ under
the inclusion
 $$
 {\left(\prod\nolimits_{i=1}^{\infty} K_0\left(\mathbb K \otimes
D_k\right)\right)/K_0\left(\mathbb K \otimes D\right)}
 \subseteq
 {\left(\prod\nolimits_{i=1}^{\infty} K_0\left(\left(\mathbb K \otimes
D_k\right)^+\right)\right)/K_0\left(\mathbb
 K \otimes D\right)}
 $$
is also the image of $P^+_*\left([e] + [f_1] - [1_n]\right) \in
{\prod_{i=1}^{\infty} K_0\left(\left(\mathbb K \otimes
D_k\right)^+\right)}$ under the quotient map
 $$
{\prod\nolimits_{i=1}^{\infty} K_0\left(\left(\mathbb K \otimes
D_k\right)^+\right)} \to  {\left(\prod\nolimits_{i=1}^{\infty}
K_0\left(\left(\mathbb K \otimes
D_k\right)^+\right)\right)/K_0\left(\mathbb K \otimes D\right) } .
 $$
Write $P^+(f_1) = \left(g_i\right)$, where $g_i \in
M_{1+n+m}\left(\left(\mathbb K \otimes D_i\right)^+\right)$ for
each $i$ is a projection whose image in $ M_{1+n+m}\left(\mathbb
C\right)$ under the canonical surjection $ M_{1+n+m}\left(
\left(\mathbb K \otimes D_i\right)^+\right) \to
M_{1+n+m}\left(\mathbb C\right)$ is a projection of rank $n$.
Since $\bigcup_j M_{1+n+m}\left( \left(M_j(\mathbb C) \otimes
D_i\right)^+\right)$ is dense in $ M_{1+n+m}\left( \left(\mathbb K
\otimes D_i\right)^+\right)$, there is a $j \in \mathbb N$ and a
projection $f^i_2 \in M_{1+n+m}\left( \left(M_j(\mathbb C) \otimes
D_i\right)^+\right)$ which is unitarily equivalent to $g_i$. Since
 $$
M_{1+n+m}\left( \left(M_j(\mathbb C) \otimes D_i\right)^+\right) =
M_{1+n+m}\left(M_j(\mathbb C) \otimes D_i\right)\oplus
M_{1+n+m}(\mathbb C),
 $$
we see that $f^i_2 = f^i_3 + f^i_4$, where $f^i_3$ and $f^i_4$ are
orthogonal projections in $ M_{1+n+m}\left( \left(\mathbb K
\otimes D_i\right)^+\right)$, $f^i_3 \in  M_{1+n+m}\left( \mathbb
K \otimes D_i\right)$, and $[f^i_4] = [1_n]$ in $
K_0\left(\left(\mathbb K \otimes D_i\right)^+\right)$. It follows
that
 $$
P^+_*\left([e] + [f_1] - [1_n]\right) = P_*[e] +
\left(\left[f^i_3\right]\right)_{i=1}^{\infty} \in
\prod\nolimits_{i=1}^{\infty} K_0\left(\mathbb K \otimes
D_k\right).
 $$
We identify now $ K_0\left(\mathbb K \otimes D_k\right)$ with
$\mathbb Z$ as ordered groups, and consequently
$\prod_{i=1}^{\infty} K_0\left({\mathbb K} \otimes D_k\right)$
with $\prod_{i=1}^{\infty} \mathbb Z$. Then $P_*\left([e]\right) =
\left(a_i\right)_{i=1}^{\infty}$ and $
\left(\left[f^i_3\right]\right)_{i=1}^{\infty} =
\left(b_i\right)_{i=1}^{\infty}$, where $b_i \geq 0$ for all $i$,
and $a_i$ is greater or equal to the multiplicity of the trivial
representation of $G$ in $\pi' \otimes \overline{\pi_i}$ which
equals $n_i$. Since $K_0\left( \mathbb K \otimes D\right)$ is the
subgroup of $\prod_{i=1}^{\infty} \mathbb Z$ consisting of the
sequences $(c_i)$ in $\mathbb Z$ for which $\sup_i
\bigl|\frac{c_i}{d_i}\bigr| < \infty$ by Lemma 3.2 of \cite{DE},
we conclude that $x \neq 0$ because $\lim_{i \to
\infty}\frac{n_i}{d_i} =\infty$. This contradicts (\ref{urt18}).

\end{proof}

 \medskip
The first named author is grateful to S. Wassermann for sending
him a copy of his book \cite{Wbook}.

\end{document}